\documentclass[12pt]{amsart}
\usepackage{amssymb}
\usepackage{amsxtra}
\usepackage[mathscr]{eucal}
\usepackage{graphics}
\usepackage{epsfig}
\usepackage{hyperref}

\newtheorem{theorem}{Theorem}[section]
\newtheorem{lemma}[theorem]{Lemma}
\newtheorem{corollary}[theorem]{Corollary}
\newtheorem{proposition}[theorem]{Proposition}

\theoremstyle{definition}

\newtheorem{definition}[theorem]{Definition}

\newtheorem{example}[theorem]{Example}

\def\PID{principal ideal domain}   %% PID
\def\Q{{\mathbb Q}}    %% rationals
\def\z{{\zeta}}        %% complex cube root of 1
\def\SE{{S_E}}         %% primes of bad reduction
\def\SK{{S_K}}         %% a finite subset of places of K

\begin{document}
\title{Rational points on elliptic curves}  
\author{Graham Everest, Jonathan Reynolds and Shaun Stevens}
\address{School of Mathematics, University of East Anglia,
Norwich NR4 7TJ, UK.}\email{g.everest@uea.ac.uk}
\dedicatory{\today}
\begin{abstract}
We consider the structure of rational points on elliptic curves in
Weierstrass form. Let~$x(P)=A_P/B_P^2$ denote the $x$-coordinate
of the rational point~$P$ then we consider when~$B_P$ can be a
prime power. Using Faltings' Theorem we show that for a fixed
power greater than~$1$, there are only finitely many rational
points with this property. 
Where descent via an isogeny is possible we show, with no
restrictions on the power, that there are only finitely many
rational points with this property, that these points are bounded in
number in an explicit fashion, and that they are effectively computable.
\end{abstract}
\thanks{The authors thank Gunther Cornelissen and Jan-Hendrik Evertse
for helpful comments.}

\maketitle

Let $E$ denote an elliptic curve given by a Weierstrass equation
\begin{equation}\label{weierstrassequation}
y^2+a_1xy+a_3y=x^3+a_2x^2+a_4x+a_6
\end{equation}
with integral coefficients $a_1,\dots, a_6$. See \cite{cassels}
and \cite{MR87g:11070} for background on elliptic curves.
Let~$E(\mathbb Q)$ denote the group of rational points on $E$. For
an element~$P\in E(\mathbb Q)$, the shape of the defining equation
(\ref{weierstrassequation}) requires that~$P$ be in the form
\begin{equation}\label{genericrationalpoint}
P=\left(\frac{A_P}{B_P^2},\frac{C_P}{B_P^3}\right)
\end{equation}
where $A_P,B_P,C_P$ are integers with no common factor. In this
paper we are concerned with the equation
\begin{equation}\label{theequation}B_P=p^f
\end{equation}
where $p$ denotes a prime and $f\ge 1$ denotes an integer. All of
the methods extend to the equation~$B_P=bp^f$, where~$b$ is fixed,
with trivial adjustments. Gunther Cornelissen once remarked to the
first author that the equation (\ref{theequation}) has only
finitely many solutions when~$f=2$, by re-arranging the equation
and invoking Faltings' Theorem.

\begin{theorem}\label{genuinetheorem}Let $E$ denote an elliptic curve
in Weierstrass form (\ref{weierstrassequation}). For any fixed~$f>1$,
there are only finitely many~$P$ for which equation
(\ref{theequation}) has a solution, with~$p$ denoting any integer.
\end{theorem}

The proof is based upon one of Siegel's proofs of his theorem
about $S$-\nobreak integral points, and uses Faltings' Theorem at a
critical stage. In his thesis \cite{jonthesis}, the second author
gives a generalization of Theorem \ref{genuinetheorem} to number
fields. 

Where descent via an isogeny is possible,
Theorem~\ref{genuinetheorem} can be strengthened. If~$G$ denotes a
subset of~$E(\mathbb Q)$, let~$r_G$ denote the rank of the
subgroup of~$E(\mathbb Q)$ generated by~$G$. Given an elliptic
curve~$E$, let~$\Delta_E$ denote the discriminant of~$E$. Finally,
let~$\omega(N)$ denotes the number of distinct prime factors of
the integer~$N$.

\begin{theorem}\label{primepower} Let $E$ denote an elliptic curve and let $G$
denote a given subset of $E(\mathbb Q)$. Assume~$G$ lies in the
image of a non-trivial isogeny from a subset~$G'$ of $E'(\mathbb Q)$, for an
elliptic curve~$E'$. There are only finitely many~$P$ for which
equation (\ref{theequation}) has a solution, where~$p$ denotes any
prime, and~$f$ denotes any integer~$f\ge 1$. The exceptional
points can be computed effectively, and they are at most
\begin{equation}\label{numberofexceptionalpoints}
\theta^{(r_{G}+1)\omega(\Delta_{E})},
\end{equation}
in number, where~$\theta
>1$ is a uniform constant which can be presented explicitly.
\end{theorem}

Note the crucial distinctions: in Theorem~\ref{primepower}
both~$p$ and~$f$ are allowed to vary whereas in
Theorem~\ref{genuinetheorem} the exponent~$f$ is fixed beforehand.
Also Theorem \ref{primepower} deals with the case~$f\ge 1$ whereas
Theorem~\ref{genuinetheorem} assumes~$f>1$. The effectiveness
statement assumes the {\em givenness} of the set~$G$;
typically~$G$ will consist of the linear span of a finite set of
points. 

The effective and explicit nature of
Theorem~\ref{primepower} renders this a much stronger outcome. The
effectiveness statement is possibly of more moral than practical
worth. The tools used include elliptic transcendence theory which,
on the face of it, exhibits a large bound for the height of the
exceptional points. In practice, the bounds are usually much
smaller. Indeed, infinite families can be constructed with very
small exceptional sets.

\begin{example}[\sc \cite{jonthesis}]\label{jonthesisexample}
For every integer~$T>1$ consider the elliptic curve
$$y^2 = (x+1)(x-T^2)(x-T^4)
$$
together with the rank-1 subgroup generated by the
point~$P=(0,T^3)$. For all~$T>1$ and all~$n>2$ the denominator
of~$x(nP)$ is divisible by at least two distinct primes. More will
be said about this example at the conclusion of Section
\ref{uniformity}.
\end{example}

Our third theorem concerns curves in homogeneous form. Suppose $E$
denotes an elliptic curve defined by an equation
\begin{equation}\label{homogform}
u^3+v^3=D,
\end{equation}
for some non-zero $D\in \mathbb Q$. Let $P$ denote a non-torsion
$\mathbb Q$-rational point. Write, in lowest terms
\begin{equation}\label{generichomogpoint}
P=\left(\frac{A_P}{B_P},\frac{C_P}{B_P}\right).
\end{equation}

\begin{theorem}\label{eveneasierfiniteness}
Suppose $E$ denotes an elliptic curve defined by an equation
(\ref{homogform}) for some non-zero cube-free $D\in \mathbb Q$. Let $P$
denote a non-torsion point in $E(\mathbb Q)$. With~$P$ as in
(\ref{generichomogpoint}), the integers $B_P$ are prime powers for
at most $\mu^{\omega(D)}$ points $P$, where~$\mu >1$ is a uniform
constant which can be presented explicitly.
\end{theorem}

\subsection{Prime Denominators}
The equation~$B_P=p$ (in other words~$f=1$) occupies something of
a middle ground between Theorems \ref{genuinetheorem} and
\ref{primepower}. The possibilities do seem to be more subtle. In
\cite{chuds}, \cite{eds}, \cite{primeds} and \cite{pe}, this
apparently more difficult question has been considered. In
\cite{primeds} and \cite{pe} we argued that, for rank-1
subgroups~$G$ (in other words, multiples of a single point) only
finitely many points $P\in G$ should yield prime values~$B_P$. In
certain cases, we could prove this and we conjectured that, for a
rank-1 subgroup~$G$, the number of rational points~$P\in G$ as in
(\ref{genericrationalpoint}), with~$B_P$ a prime is uniformly
bounded if the defining equation (\ref{weierstrassequation}) is in
minimal form. Many calculations suggest that for an elliptic curve
in minimal form, a little over $30$ prime values $B_P$ seems to be
the limit. The following example may well provide the limiting case.
It is taken from Elkies' table of small height points in
\cite{elkies}. Note that the curves in Elkies' table are not
presented in minimal form; the example following is the second on
the list, with the equation rendered in minimal form. The point
shown has the second smallest known height amongst all rational
points on elliptic curves.
\begin{example}Let $P$ denote the point $P=(-386,-3767)$ on the
elliptic curve
$$
y^2+xy=x^3-141875x+13893057.
$$
The values $x(nP)$ have prime square denominators for~$31$ values
of~$n$, ending with (apparently)~$n=613$. The first example on
Elkies' list yields (apparently)~$30$ such primes. These computations
were performed using PARI-GP.%~\cite{parigp}.
\end{example}
In a short subsection (\S\ref{uniformity}) we will give an outline of a
proof of uniformity assuming an unproven (yet generally believed
conjecture) as well as an improvement to known results in elliptic
transcendence theory.

On the other hand, in \cite{rew}, many examples of rank-2 curves
were considered where, apparently, there are infinitely many
rational points~$P$ with~$B_P$ prime.
\begin{example}\label{with}The curve
$$y^2=x^3-28x+52
$$
has rank 2, with generators $P_1=(-2,10)$ and $P_2=(-4,10)$. In
\cite{rew} we considered the possibility that asymptotically $\rho
\log T$ values $$x(n_1P_1+n_2P_2) \mbox{ with } \max
\{|n_1|,|n_2|\}<T$$ possess a prime square denominator, where
$\rho
>0$ is a constant depending upon $E$.
The table in \cite{2deds}, constructed by Peter Rogers, exhibits
many more rank-2 curves which are thought to yield infinitely many
rational points~$P$ as in (\ref{genericrationalpoint}) with~$B_P$
a prime.
\end{example}
Theorem~\ref{primepower} allows many examples which yield only
finitely many primes to be constructed, when a descent is possible.

\begin{example}\label{morewithout}The curve
$$
y^2=x^3-6400x+192000
$$
has rank 2, with independent points $P_1=(65,225)$ and
$P_2=(56,96)$. Theorem \ref{primepower} guarantees there are
finitely many prime square denominators amongst the points
$n_1P_1+n_2P_2$ because the points shown are the images of the
points $(0,100)$ and $(6,104)$ under a 2-isogeny from the curve
$$
y^2 = x^3+100x+10000.
$$
\end{example}

In Section \ref{primepowers} we give a proof of Theorem
\ref{genuinetheorem}: since it relies on Faltings' Theorem, it is
not effective. Following this, we prove Theorem \ref{primepower};
the effectiveness statement uses local heights and elliptic
transcendence theory whilst the explicit bound for the exceptional
points relies upon a strong version of Siegel's Theorem
\cite{grosssilverman}, \cite{silvermansiegel}. The final section
proves Theorem~\ref{eveneasierfiniteness}, using a strong form of
Thue's Theorem to rid us of the dependence upon the rank.

\section{Proof of Theorem \ref{genuinetheorem}}\label{primepowers}

%\begin{proof}[\sc Proof of Theorem (\ref{genuinetheorem})]
Let $K$ be a finite field extension of $\mathbb{Q}$. The
valuations on $K$, written~$M_K$, consist of the usual archimedean
absolute values together with the non-archimedean, $\wp$-adic
valuations, one for each prime ideal $\wp$ of~$K$. Write $K_v$ for
the completion of $K$ with respect to $v$. Let $S \subset M_K$
denote a finite set of valuations containing the archimedean
valuations. The ring $O_S$ of $S$-integers is given by
\[
O_S = \{x\in K : \nu(x) \ge 0 \textrm{ for all } \nu \in M_K, \nu
\notin S \}
\]
and the unit group $O_S^*$ of $O_S$ is given by
\[
O_S^* = \{x\in K : \nu(x) = 0 \textrm{ for all } \nu \in M_K, \nu
\notin S \}.
\]
Completing the square in (\ref{weierstrassequation}), it is
sufficient to consider an equation
\begin{equation} \label{nop}
y^2=x^3+a_2x^2+a_4x+a_6
\end{equation}
where $x^3+a_2x^2+a_4x+a_6 \in \mathbb{Q}[x]$ has distinct zeros 
$\alpha_1, \alpha_2, \alpha_3$ in some finite extension $K$ of
$\mathbb Q$. Of course we might have introduced powers of $2$ into the
denominators in \eqref{nop}; we will see that this does not matter. 

Let $S$ be a sufficiently large (finite) subset of $M_K$ so that $O_S$
is a \PID and $2,\alpha_i-\alpha_j \in O_S^*$ for all $i \ne j$. 
Now let $L/K$ be the extension of $K$ obtained by adjoining to $K$
the square root of every element of $O_S^*$. Note that $L/K$ is a
finite extension, since $O_S^*/(O_S^*)^2$ is finite from
Dirichlet's $S$-unit theorem. Further let $T \subset M_L$ be a
finite set containing the places of $L$ lying over elements of
$S$ and such that $O_T$ is a \PID, where, by abuse of notation, $O_T$
denotes the ring of $T$-integers in $L$.

Now we turn to the proof of the Theorem. Replacing $(x,y)$ by
$(x/q^2,y/q^3)$ in \eqref{nop}, we are searching for 
solutions in $\Q\cap O_S$ to
\begin{equation} \label{p}
y^2=x^3+a_2q^2x^2+a_4q^4x+a_6q^6,\qquad \gcd(xy,q)=1.
\end{equation}
We will show that, for fixed $f>1$, there are only finitely many
prime powers $q=p^f\in \mathbb{Z}$ for which \eqref{p} has a solution. Note
that, since $f$ is fixed and $T$ is finite, we may assume that $p$ is
large enough so that no valuation of $L$ dividing $p$ lies in $T$.

Let $x,y\in \mathbb Q\cap O_S$ be a
solution to \eqref{p}; then
\begin{equation} \label{sq}
y^2=(x-q^2\alpha_1)(x-q^2\alpha_2)(x-q^2\alpha_3).
\end{equation}
Let $\wp$ be a
prime ideal of $O_S$ dividing $x-q^2\alpha_i$; then $\wp$ cannot
divide $q$, since $(x,q)=1$. Hence $\wp$ can divide at most one term
$x-q^2\alpha_i$, since if it divides both $x-\alpha_iq^2$ and
$x-\alpha_jq^2$ then it divides also $(\alpha_i-\alpha_j)q^2$. From
\eqref{sq} it follows that there are elements $z_i \in O_S$ and
units $b_i \in O_S^*$ so that
\[
x-\alpha_iq^2=b_iz_i^2.
\]
We have $b_i= \beta_i^2$, for some $\beta_i\in O_T$ so
\begin{equation} \label{3}
x-\alpha_iq^2=(\beta_iz_i)^2.
\end{equation}
Taking the difference of any two of these equations yields
\[
(\alpha_j-\alpha_i)q^2=(\beta_iz_i-\beta_jz_j)(\beta_iz_i+\beta_jz_j).
\]
Note that $\alpha_j-\alpha_i \in O_T^*$ while each of the two
factors on the right is in $O_T$. It follows that each of these
factors is made from primes $\pi | p$ in $O_T$. Further we may
assume these factors are coprime, since if $\pi| p$ divides
$2\beta_iz_i$ and $p>2$ then from (\ref{3}) $\pi$ divides $x$. 
%If
%$\beta_iz_i+\beta_jz_j$ is a unit for all $i \ne j$ then Siegel's
%identity:
%\[
%\frac{\beta_1z_1+\beta_2z_2}{\beta_1z_1-\beta_3z_3}-\frac{\beta_2z_2+\beta_3z_3}{\beta_1z_1-\beta_3z_3}=1
%\]
%gives
%\[
%u+v=q^2,
%\]
%for $S$-units $u$ and $v$. 

Without loss of generality, we assume that
there is a prime $\pi \in O_T$ dividing $\beta_1z_1+\beta_2z_2$.
If $\pi$ does not divide $\beta_2z_2+\beta_3z_3$ then $\pi$ does
not divide $\beta_1z_1-\beta_3z_3$. So Siegel's identity 
\[
\frac{\beta_1z_1+\beta_2z_2}{\beta_1z_1-\beta_3z_3}-
\frac{\beta_2z_2+\beta_3z_3}{\beta_1z_1-\beta_3z_3}=1
\]
gives
\begin{equation} \label{4}
a_p^2u+b_p^2v=c_p^2
\end{equation}
where $a_p,b_p,c_p \in O_T$ divide $q$, are not all $T$-units and
are pairwise coprime. If $\pi$ divides $\beta_2z_2+\beta_3z_3$
then $\pi$ divides $\beta_1z_1-\beta_3z_3$ so $\pi$ does not
divide $\beta_1z_1+\beta_3z_3$. Then Siegel's identity
\[
\frac{\beta_2z_2+\beta_1z_1}{\beta_2z_2-\beta_3z_3}-\frac{\beta_1z_1+\beta_3z_3}{\beta_2z_2-\beta_3z_3}=1
\]
gives (\ref{4}).

Since $q=p^f$ is a prime power with~$f>1$, each of $a_p$ and
$b_p$ is itself an $f^{\rm th}$ power. Finally, we note that the
group~$O_T^*/(O_T^*)^{2f}$ is finite so we fix once and for 
all a set of coset representatives. Then, by \eqref{4}, each solution
to \eqref{p} gives us a solution of an equation:
$$
ux^{2f}+vy^{2f}=1, \qquad x,y\in L,
$$
with $2f \ge 4$, where~$u$ and~$v$ belong to this finite set of
representatives, which
depends only upon~$T$ and~$f$. Since each such curve has genus
$(2f-1)(f-1)\ge 3$, Faltings' Theorem \cite{ft} guarantees there are
only finitely many solutions. Since we have only finitely many such
equations, we are done.
\hfill$\square$
%\end{proof}

\section{Proof of Theorem \ref{primepower}}

Let $E$ and $E'$ be two elliptic curves, defined over
$\mathbb Q$. An {\em isogeny} is a non-zero homomorphism
\begin{equation*}
\sigma :E'\rightarrow E
\end{equation*}
taking the zero of $E'$ to the zero of $E$. There is a dual
isogeny $\sigma^*:E\rightarrow E'$ and the composite homomorphisms
$\sigma\sigma^*$ and $\sigma^*\sigma$ are multiplication by $d$ on $E$
and $E'$ respectively, for some integer $d$, which is said to be
the \textit{degree} of the isogeny. The curves $E$ and $E^{\prime
}$ are said to be {\em $d$-isogenous} if there is an isogeny of
degree $d$ between them.

\begin{definition}\label{defofS}For any
point $P\in E(\mathbb Q)$ let $S(P)$ denote the non-archimedean
valuations $v$ in $M_{\mathbb Q}$ for which $|x(P)|_v>1$. The
definition of $S(P)$ depends upon the Weierstrass equation so we
assume this equation has been fixed.
\end{definition}

\begin{theorem}\label{finite} Let $E$ denote an elliptic curve
which is defined over $\mathbb Q$ and let $G \subset E(\mathbb Q)$
denote a subset contained in the image of a subset~$G'$ of
rational points under a non-trivial isogeny. Then there are only finitely
many $P \in G$ for which $S(P)$ consists of a single valuation:
these points are effectively computable and they are at most
$c^{(r_{G}+1)\omega(\Delta_{E})}$ in number.
\end{theorem}

Theorem \ref{finite} will be proved in Section \ref{proof}
following a section with some basic properties of heights under
isogeny.

\subsection{Heights}\label{hival}

Write
\begin{equation}\label{localheight}
h_v(\alpha)=\log \max \{ 1, |\alpha|_v\},
\end{equation}
for the local (logarithmic) height at $v$. The na\" ive global
logarithmic height of $Q$ is defined to be
$$h(\alpha)=\sum_{v\in M_{\mathbb Q}} h_v(\alpha)=
\sum_{v\in M_{\mathbb Q}} \log \max \{ 1, |\alpha|_v\},
$$
the sum running over all the valuations of $\mathbb Q$. If $P$
denotes any rational point on an elliptic curve, we write 
$h_v(P)=h_v(x(P))$ and $h(P)=h(x(P))$. Usually, in the
literature, the global height is further normalized by dividing by
$2$.

Suppose $P$ denotes a rational point of $E$. The theory of heights
gives an estimate for
$$
h(P)=\widehat h(P)+O(1),
$$
where $\widehat h(P)$ denotes the canonical height of $P$. The
canonical height enjoys the additional property that $\widehat
h(mP)=m^2\widehat h(P)$ for any $m\in \mathbb Z$. More generally,
if $\sigma:E'\rightarrow E$ denotes a $d$-isogeny then for $P'\in
E'(\mathbb Q)$
\begin{equation}\label{goingup}\widehat h(\sigma(P'))=d\widehat h(P').
\end{equation}

\begin{lemma}\label{growth}
Suppose $P_1,\dots,  P_r \in E(\mathbb Q)$ are independent
rational points and $T$ is a torsion point. Then
$$
h_v(n_1P_1+\dots +n_rP_r+T)=
O(\log |\underline n| (\log \log |\underline n|)^{r+2}),
$$
for any valuation $v$\ where $|\underline n|=\max \{|n_1|,\dots
,|n_r|\}$ for $\underline n \in \mathbb Z^r$. This can
be written
\begin{equation}\label{growthagain}
\log |x(P)|_v =O(%\ll 
\log \widehat h(P) (\log \log \widehat
h(P))^{r+2}),
\end{equation}
for any $P\in E(\mathbb Q)$.
\end{lemma}

\begin{corollary}\label{fixedset}
With the same notation as Lemma~\ref{growth}, let
$S$ be any fixed, finite set of valuations of
$\mathbb Q$. Then
$$
\sum_{v\in S}h_v(n_1P_1+\dots +n_rP_r+T)=O((\log |\underline n|)^2).
$$
This can be written
\begin{equation}\label{finitegrowth}
\sum_{v\in S}h_v(P)=O((\log \widehat h(P))^2),
\end{equation}
for any $P\in E(\mathbb Q)$.
\end{corollary}

\begin{proof}[Proof of Lemma~\ref{growth}.]
We will detail a proof for the archimedean valuation only. The
proof for non-archimedean valuations is similar. The estimate in
Lemma~\ref{growth} follows from an appropriate upper bound for
$|x(n_1P_1+\dots +n_rP_r)|_v$. Putting the given model of the
curve into a short Weierstrass equation only translates $x$ by at
most a constant. Let $z_{P_i}$ correspond to $P_i$ under an
analytic isomorphism $E(\mathbb C)\simeq \mathbb C/L$, for some
lattice $L$, with~$z_T$ corresponding to~$T$. Thus we may assume
that the $x$-coordinate of a point is given using the Weierstrass
$\wp$-function with Laurent expansion in even powers of $z$,
$$x=\wp_L(z)=\frac{1}{z^2}+c_0+c_2z^2+\dots
$$
Write $\{n_1z_{P_1}+\dots +n_rz_{P_r}+z_T\}$ for $n_1z_{P_1}+\dots
+n_rz_{P_r}+z_T$ modulo $L$. When the quantity $|x(n_1P_1+\dots
+n_rP_r+T)|$ is large it means $n_1P_1+\dots +n_rP_r+T$ is close
to zero modulo $L$, thus the quantities $|x(n_1P_1+\dots
+n_rP_r+T)|$ and $1/|\{n_1z_{P_1}+\dots +n_rz_{P_r}+z_T\}|^2$ are
commensurate. On the complex torus, this means the elliptic
logarithm is close to zero. So it is sufficient to supply a lower
bound for $|\{n_1z_{P_1}+\dots +n_rz_{P_r}+z_T\}|$ and this can be
given by elliptic transcendence theory (see \cite{dav}). We use
Th\'eor\`eme 2.1 in~\cite{dav} but see also~\cite{ST} where an
explicit version of David's Theorem appears on page 20. The nature
of the bound is
\begin{equation}\label{b2}
\log |x(n_1P_1+\dots +n_rP_r+T)| \ll \log |\underline n| (\log
\log |\underline n|)^{r+2} ,
\end{equation}
where the implied constant depends upon $E$, the valuation $v$ and
the points $P_1,\dots ,P_r$.

For the final assertion, we need only note that the global canonical height 
$\hat h(n_1P_1+\dots +n_rP_r+T)$ is a positive definite quadratic form
in $n_1,...,n_r$, and hence comensurate with $|\underline n|^2$.
\end{proof}

We will need some more theory of elliptic curves over local
fields, see~\cite{MR87g:11070}. For every non-archimedean
valuation $v$, write ord$_v$ for the corresponding order function.
There is a subgroup of the group of $\mathbb Q_{v}$-rational
points:
$$
E_1(\mathbb Q_{v})=\{O\} \cup \{P\in E(\mathbb Q_{v}):\mbox{ord} _{v}(x(P))
\leq -2\}.
$$ 
In \cite{MR87g:11070}, Silverman proves the
following.

\begin{proposition}For all $P\in E_1(\mathbb Q_{v})$ and all $d\in \mathbb Z$:
\begin{equation}\label{ordthing}
\log|x(mP)|_v=\log |x(P)|_v-\log |m|_v.
\end{equation}
\end{proposition}

For finitely many (bad) primes $p$, the reduction of $E$ is not an elliptic
curve because the reduced curve is singular. We write $S_E$ for the
set of valuations in $M_\Q$ corresponding to all such primes. The
equation~\eqref{ordthing} then yields the following corollary:

\begin{corollary}\label{bigOthing} Suppose $\sigma:E'\rightarrow E$
is a $d$-isogeny and $P'\in E'_1(\mathbb Q_v)$. If $v \notin \SE$
then $h_v(\sigma(P'))\ge h_v(P')$. The local heights are related
by the formula
\begin{equation}\label{newordthing}
h_v(\sigma(P'))=h_v(P')+O(1),
\end{equation}
where the implied constant depends only upon the isogeny and is
independent of $P'$.
\end{corollary}

\begin{proof} %[Proof of Corollary \ref{bigOthing}]
Suppose $v$ corresponds to the prime $p$. Provided $v \notin \SE$,
both curves and the isogeny can be reduced modulo powers of $p$
and the first statement in the Corollary follows. Applying the
dual isogeny $\sigma^*$ gives a similar inequality $h_v(\sigma^*(P))\ge
h_v(P)$, for all $P\in E_1(\mathbb Q_v)$. However, composing
$\sigma$ with its dual gives multiplication by $d$ on $E'$. Now
(\ref{ordthing}) applies to prove (\ref{newordthing}).
\end{proof}

\subsection{Proof of Theorem \ref{finite}}\label{proof}

\hfill\break
%\begin{proof}[Proof of Theorem \ref{finite}]
Suppose $G'$ is a subset of $E'(\mathbb Q)$ and $\sigma:E'\rightarrow
E$ is a $d$-isogeny with $\sigma(G')=G$. From
(\ref{goingup}),
\begin{equation}\label{compareheights}
\widehat h(P)=d\widehat h(P'),
\end{equation}
where $\widehat h$ denotes the canonical heights of rational
points on~$E$ and~$E'$. 

Suppose $S(P)$ consists of the single valuation $v$. 
If $v\in\SE$ then, by~\eqref{finitegrowth}, we have 
$$
\widehat h(P')=O((\log \widehat h(P'))^2).
$$ 
If $v\not\in\SE$ then $P\in E_1(\Q_v)$ and, by reduction,
$S(P')\subset S(P)$. Again by~\eqref{finitegrowth},
\begin{equation}\label{heightdifference}
h(P)-h_v(P)=O((\log \widehat h(P'))^2)=h(P')-h_v(P').
\end{equation}
Now the canonical height differs from the na\"ive height by a bounded
amount so we are justified in using the canonical height
in~\eqref{heightdifference}. It follows from~\eqref{compareheights}
that
$$
d\widehat h(P')-h_v(P)=O((\log \widehat h(P'))^2).
$$
However, from~\eqref{newordthing}
$$
h_v(P')-h_v(P)=O(1).
$$
Subtracting these last two formulae, and dividing by $d-1>0$, gives
\begin{equation}\label{heightlessthanlogheight}
\widehat h(P')=O((\log \widehat h(P'))^2),
\end{equation}
with an implied constant which is effectively computable. In
particular, we have obtained the same inequality for any valuation
$v$. Equation~\eqref{heightlessthanlogheight} bounds the height
$\widehat h(P')$ so can only be satisfied 
by finitely many points~$P'$, which can be computed effectively.

\medskip

Finally, we estimate the number of exceptional points, using a
strong form of Siegel's Theorem, proved by Gross and Silverman.
There is a non-trivial torsion point~$T'$ in the kernel
of~$\sigma:E'\rightarrow E$. Any isogeny factorizes as a product
of isogenies of prime degree so, since the rank $r_G$ and the
discriminant $\Delta_E$ are preserved under isogeny, we may assume from
now that~$\sigma$ has prime degree. Let~$K$ denote a number field over
which~$T'$ is 
defined. By Mazur's famous result (\cite[Theorem 1]{mazur}), only
finitely many primes can occur as degrees of prime degree isogenies
which map onto rational points (the largest of which is
$163$). Let~$S$ denote the subset of $M_\Q$ consisting of all such
prime degrees, together with the primes of bad reduction (in $\SE$). 
Let $\SK$ denote the subset of $M_K$ consisting of all places above
those in $S$, together with the places dividing the denominator of
$T'$. 

Now the points~$P'$ and~$P'+T'$ both map to~$P$ under~$\sigma$. It
follows that if ~$P'$ (respectively ~$P'+T'$) has denominator divisible
by a place of good reduction $w$ (respectively $w'$), then $w$
(respectively $w'$) is guaranteed to
appear in the denominator of~$P$. Moreover, if~$w,w'\not\in\SK$, they
are guaranteed to be distinct: indeed, any
good reduction place dividing the denominator of both
points will divide the denominator of~$T'$ and hence lie in
$\SK$. Further, since $P'$ is a $\Q$-rational
point, $w'$ is coprime to the prime of $\Q$ below $w$ so there are two
distinct primes of $\Q$ dividing~$B_P$.

This ensures~$B_P$ cannot be a prime power unless either~$P'$
or~$P'+T'$ is an~$\SK$-integral point on~$E'$. By the
theorem of Gross and Silverman \cite[Theorem 0.1]{grosssilverman},
there are at most 
$$
d\eta^{r\delta (j_{E'})+|\SK|}
$$
$\SK$-integral points in~$E'(K)$, for an explicit constant~$\eta$,
where~$d=[K:\mathbb Q]$ and~$\delta (j_{E'})$ denotes the number
of primes dividing the denominator of the $j$-invariant of~$E'$.
In their paper,~$r$ denotes the rank of the group~$E(K)$; however,
their results are valid for the number of $\SK$-integral points
inside a subgroup of rank~$r$. Since the primes dividing the
denominator of~$j_{E'}$ divide~$\Delta_{E'}=\Delta_E$, since~$|\SK|$ is a 
bounded multiple of~$d\omega(\Delta_{E})$, and since~$d$ is
uniformly bounded by Mazur's Theorem, the bound in
(\ref{numberofexceptionalpoints}) follows.
\hfill$\square$
%\end{proof}

\subsection{Uniformity}\label{uniformity}

Under the assumption that $G$ consists of a rank-1 subgroup
of~$E(\mathbb Q)$, we may harness the ideas in Section \ref{proof}
to explain how these might be used as part of a proof of
uniformity where a descent is possible. In the rank-1 case, the
inequality (\ref{heightlessthanlogheight}) can be stated more
explicitly. Let~$P$ denote a generator of~$G$ and~$\sigma(P')=P$.
Writing~$h'=\widehat h(P')$ and invoking the explicit form of
David's Theorem in \cite{ST} shows that the integers~$n$ for
which~$S(nP)$ consists of a single valuation must satisfy
\begin{equation}\label{beexplicit}h'n^2<\tau \log n(\log \log n+\log
\Delta_{E'})^3,
\end{equation}
where $\tau$ denotes a uniform bound. Lang's Conjecture asserts a
uniform upper bound for~$(\log \Delta_{E'})/h'$. If the dependence
upon~$\log \Delta_{E'}$ on the right hand side of
(\ref{beexplicit}) were linear, then Lang's Conjecture would
guarantee a uniform upper bound for the number of prime square
denominators in the sequence~$x(nP)$ when descent is possible.

The reason Example \ref{jonthesisexample} works is that Lang's
Conjecture can be proved in an explicit manner for 1-parameter
families such as this \cite{silvariation}. Also, it is possible to
obtain a stronger form of Lemma \ref{growth} without transcendence
theory: since~$P$ lies off the connected component of the identity
so do all its odd multiples and for these an upper bound for the
$x$-coordinate exists. Also, (\ref{ordthing}) allows an easy
treatment of the bad reduction primes by showing the local
heights~$h_v(nP)$ are each~$O(\log n)$.

\section{The curve $u^3+v^3=D$}\label{othermodels}

In \cite{primeds}, it was proved that the integers $B_P$ are prime
powers for only finitely many $\mathbb Q$-points $P$. The proof
used the well-known bi-rational equivalence of (\ref{homogform})
with the curve
$$
y^2=x^3-432D^2.
$$

{\bf Example} As Ramanujan famously pointed out, the taxi-cab
equation
\begin{equation}\label{taxicab}
u^3+v^3=1729,
\end{equation}
has two distinct integral solutions. These give rise to points
$P=[1,12]$ and $Q=[9,10]$ on the elliptic curve (\ref{taxicab}).
The only rational points on (\ref{taxicab}) which seem to yield
prime power denominators are $2Q$ and $P+Q$ (and their inverses).

\begin{proof}[Proof of Theorem \ref{eveneasierfiniteness}]
We assume that $D$ is integral; if $D$ is not integral then we can
scale by a non-zero integer to reduce to this case. 
Firstly, we recall the easy proof of Siegel's Theorem for curves
in homogeneous form~(\ref{homogform}). The equation factorizes as
$$
(u + v)(u^2 - uv + v^2) = D
$$ 
so 
$$ 
|u^2 - uv + v^2| \leq |D|.
$$
But the left hand side of this is $(u - v/2)^2 + 3v^2/4$ so $|v|^2
\le 4|D|/3$. Hence the result with an explicit bound for $|v|$.
The bound for $|u|$ is identical.

This bounds the number of solutions as $O(|D|^{1/2})$ but we can
easily improve 
on this. First we have $u+v=m$, and $u^2-uv+v^2=n$, for some
integers $m,n$ such that $mn=D$ and $\gcd(m,n)|3$. 
If $\z$ is a non-trivial cube root of unity, then we get
$$
(u+\z v)(u-\z v)= n,
$$
where the factors on the left hand side
have greatest common divisor dividing $2$. The number of ways of
factorizing $n$ in this way is $O(2^{2\omega(n)})$ so the total number
of solution $(u,v)$ is $O(\mu^{\omega(D)})$, for some (explicit)
uniform constant $\mu$.

\medskip

Now suppose
$$ 
(u + v)(u^2 - uv + v^2) = u^3 + v^3 = Dq^3
$$
where $q=p^f$ is a prime power and $u,v$ are integers coprime to
$q$. We consider first the case $p>3$.
Then $q$ cannot possess a factor in common with both
brackets. Hence one bracket is $m$ and the other is $nq^3$, where
$mn=D$ and $\gcd(m,n)|3$. If the quadratic bracket
is bounded then we bound the number of solutions as before, so we
assume
$$ 
u + v = m \mbox{ and } (u+\z v)(u-\z v) = nq^3.
$$
%If $p$ is inert in $\Q(\z)$ then $q$ cannot possess a factor in common
%with both $u+\z v$ and $u-\z v$ so one of them divides $n$ in $\Z[\z]$
%and we bound the number of solutions again as before.
%So we assume $p$ is split in $\mathbb Q (\z)$, in which case we are
%looking at 
Since $u+\z v$ and $u-\z v$ are coprime, we get
$$
u + v = m \mbox{ and } u - \z v = k\rho^3
$$
where $k$ divides $n$ and $\rho$ divides $q$ in $\mathbb Z[\z]$. If we write $k
= c + d\z$ and $\rho = a + b\z$, with $a,b,c,d \in \mathbb Z$, then we
can write $u$ and $v$ explicitly in terms of $a,b,c,d$. 
Substituting into the equation $u + v = m$, we get
$$ 
(c + d)a^3 + (3c - 6d)a^2b + (3d - 6c)ab^2 + (c + d)b^3 = m.
$$
For each of the finitely many values of $c$ and $d$, this is a
Thue Equation and it is non-singular -- that is, for any
non-zero $c$ and $d$, the cubic
$$(c + d)X^3 + (3c - 6d)X^2 + (3d - 6c)X + (c + d)
$$
does not have repeated roots, since its discriminant is
$$
729c^4 - 1458dc^3 + 2187d^2c^2 - 1458 d^3c + 729d^4 =729(c^2 - cd
+ d^2)^2,
$$ 
which does not vanish unless $c = d = 0$. Thus each of
the finitely many Thue equations has finitely many solutions so
there were only finitely many values of $q$.

By~\cite[Theorem 1]{ev} (see also \emph{op. cit.} page 122), a
non-singular integral cubic Thue Equation
$$
F(x,y)=m,
$$
with $m$ cube-free,
has a number of integral solutions which is bounded by
$\mu^{\omega(m)}$, for some (explicit) uniform constant $\mu$. This
must be multiplied by the total number of 
equations, which depends only upon the number of factorizations of $mn=D$
with $\gcd(m,n)|3$, so does not change the shape of the bound claimed
by the theorem.

Finally, the cases $p=2,3$ are dealt with similarly, with minor
alterations. The point is that the greatest common divisors of the
various brackets are uniformly bounded, so the shape of the final number
of solutions is unchanged.
\end{proof}

%\bibliographystyle{plain}
%\bibliography{../bib/refs}

\end{document}